# Notes on Beltrami's *Essay*

**Steven Rose, University College London, May 2026**

## 1. The Metric of a Pseudospherical Surface

In his *Essay on the Interpretation of Non-Euclidean Geometry*, published in 1868, Eugenio Beltrami (1835-1900) stated that surfaces of constant negative curvature (which he referred to generically as pseudospheres) can be mapped onto a disc in the plane, geodesics on the pseudosphere appearing as segments of chords on the disc. More importantly he showed that the identities of hyperbolic geometry satisfied by figures on the pseudosphere are also satisfied by figures on the disc.

He began by stating without proof the metric of a surface of negative curvature

$$ds^2 = R^2 \left[ (a^2 - v^2)\, du^2 + 2uv\, du\, dv + (a^2 - u^2)\, dv^2 \right] / (a^2 - u^2 - v^2)^2 , \qquad (1)$$

where the curvature of the pseudospherical surface is equal to $-1/R^2$. This metric can be expressed in terms of Gauss' first fundamental form as

$$ds^2 = E\, du^2 + 2F\, du\, dv + G\, dv^2 .$$

He noted that if $U, V$ are the parametrised coordinates of the pseudosphere while $X, Y$ denote the coordinates of the auxiliary plane, the equations

$$U = X$$

$$V = Y$$

establish a bijection in which the pseudospherical surface is mapped onto the inside of a circle of radius $a$ whose centre lies at the origin.

The metric given in (1) above can be derived using projective geometry, introduced into this field by Felix Klein writing a few years after Beltrami (see Appendix 1).

## 2. The Metric in Euclidean Polar Coordinates

Beltrami also expressed this metric in the form

$$\boldsymbol{ds^2 = R\,[w^2(du^2 + dv^2) + (u\,du + v\,dv)^2\,]\,/\,w^2\,,} \qquad \dots\ (2)$$

where $w^2 = (a^2 - u^2 - v^2) = (a^2 - r^2)$.

This form can be used to derive the metric in Euclidean polar coordinates $r, \phi$:

Let $u = r\cos\phi$, $v = r\sin\phi$.

Then

$$du = \cos\phi\,dr - r\sin\phi\,d\phi$$

$$dv = \sin\phi\,dr + r\cos\phi\,d\phi.$$

Therefore

$$du^2 + dv^2 = dr^2 + r^2\phi^2$$

and

$$(u\,du + v\,dv)^2 = r^2 + r^2dr^2,$$

Hence

$$ds^2 = R^2\,[(dr^2 + r^2\phi^2)\,/\,(a^2 - r^2) + r^2dr^2\,/\,(a^2 - r^2)^2]$$

or

$$\boldsymbol{ds^2 = R^2\,\{a^2\,dr^2\,/\,(a^2 - r^2)^2 + r^2\,d\phi^2\,/\,(a^2 - r^2)\,\}} \qquad \dots\ (3)$$

## 3. The Angle between Geodesic Curves

Beltrami noted that the angle $\theta$ between the coordinate curves $u$ and $v$ is given by the formulas

$$\boldsymbol{\cos\theta = uv\,/\,[(a^2 - u^2)\,(a^2 - v^2)]^{1/2}} \qquad \dots\ (4)$$

$$\boldsymbol{\sin\theta = a\,(a^2 - u - v^2)^{1/2}/\,[(a^2 - u^2)\,(a^2 - v^2)]^{1/2}} \qquad \dots\ (5)$$

These formulas can be derived from the metric as follows. The coefficients $E$, $F$ and $G$ of the first fundamental form are

$$E = R^2/\,(a^2 - v^2)\,/\,(a^2 - u^2 - v^2)^2$$

$$F = R^2\,uv\,/\,(a^2 - u^2 - v^2)^2$$

$$G = R^2/\,(a^2 - u^2)\,/\,(a^2 - u^2 - v^2)^2.$$

The angle between coordinate curves $u$ and $v$ is given in terms of their cosine and sine by the well known formulas

$$\cos\theta = F / \sqrt{(EG)}$$

$$= uv / \sqrt{[(a^2 - v^2)(a^2 - u^2)]}$$

and

$$\sin\theta = [(EG - F^2) / (EG)]^{1/2}$$

$$= a (a^2 - u - v^2)^{1/2} / [(a^2 - u^2)(a^2 - v^2)]^{1/2}.$$

Beltrami noted that if either $u = 0$ or $v = 0$, equation (4) implies that the angle $\theta$ between them will be 90°. So any geodesic $u$ = constant will be orthogonal on the disc to the diameter $v = 0$ and vice-versa.

## 4. The Relation between Euclidean and Hyperbolic Distance on the Disc

For the arc of a geodesic containing the point $u = v = 0$, Beltrami gave the relation between its hyperbolic length $\rho$ and its Euclidean length $r$ as

$$\boldsymbol{\rho = R/2 \log [(a + r) / (a - r)]} \qquad \textbf{... (6)}$$

He derived this identity from the metric given in equation (1):

$$ds^2 = R^2 [(a^2 - v^2)\, du^2 + 2uv\, du\, dv + (a^2 - u^2)\, dv^2] / (a^2 - u^2 - v^2)^2.$$

Assuming that

$$u = r \cos\mu$$

$$v = r \sin\mu,$$

it follows when $\mu$ is constant that

$$du = \cos\mu\, dr$$

$$dv = \sin\mu\, dr .$$

Writing $ds$ as $d\rho$ and $r^2$ as $(u^2 + v^2)$, equation (1) then becomes

$$d\rho^2 = R^2 [(a^2 \cos^2\mu\, dr^2 - r^2 \sin^2\mu \cos^2\mu\, dr^2 + 2r^2 \sin^2\mu \cos^2\mu\, dr^2 + a^2 \sin^2\mu\, dr^2 - r^2 \sin^2\mu \cos^2\mu\, dr^2] / (a^2 - r^2)^2$$

$$= R^2 a^2 / (a^2 - r^2)^2,$$

whence

$$d\rho = R\, a / (a^2 - r^2).$$

Integration then gives

$$\rho = R/2 \, log \, [(a + r) / (a - r)],$$

as required. This formula ensures that the hyperbolic distance $\rho$ of a point from the origin becomes infinitely great as its Euclidean distance $r$ approaches $a$, the radius of the disc, allowing the whole of the hyperbolic plane to be mapped to the disc.

The distance formula (6) can also be written as

$$\boldsymbol{\rho = R \, tanh^{-1} (r/a)} \qquad \textbf{... (7)}$$

which implies, as Beltrami noted, that

$$\boldsymbol{r = a \, tanh \, (\rho / R).} \qquad \textbf{... (8)}$$

A further identity presented by Beltrami is

$$\boldsymbol{cosh \, (\rho / R) = a / (a^2 - r^2)^{1/2}} \qquad \textbf{... (9)}$$

This identity can be derived as follows. Dividing both sides of equation (6) gives

$$\rho / R = ½ \, log \, [(a + r) / (a - r)]$$

$$= log \sqrt{[(a + r) / (a - r)]}.$$

Hence

$$cosh \, (\rho / R) = [e^{log \sqrt{[(a+r)/(a-r)]}} + e^{-log \sqrt{[(a+r)/(a-r)]}}] / 2$$

$$= ½ \sqrt{[(a + r) / (a - r)]} + ½ \sqrt{[(a - r) / (a + r)]}$$

$$= ½ \, [(a + r) + (a - r)] / (a^2 - r^2)^{1/2}$$

$$= a / (a^2 - r^2)^{1/2}.$$

**5. Geometrical Derivation of the Formula for Hyperbolic Distance**

The formula for the relation between hyperbolic and Euclidean distance on the disc given in equation (7) can also be derived geometrically using Cayley's metric, based on the cross ratio. Given four collinear points on a disc centred at O with radius $a$, P, O, B, Q, where P, Q are the extremities of the diameter on which OB lies,

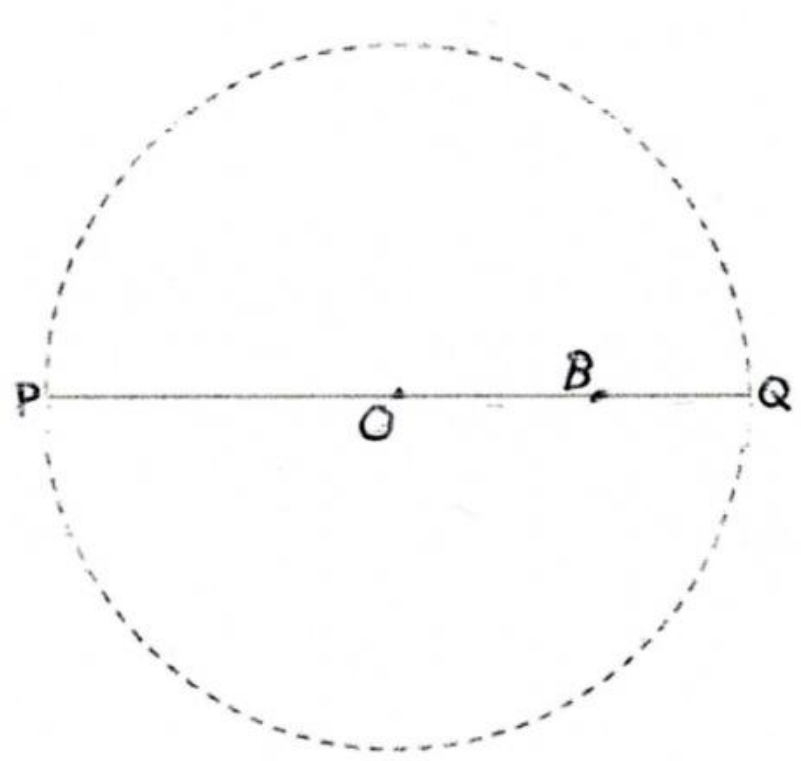

the cross ratio (OBPQ) is given by

$$(OBPQ) = PB \,.\, OQ \,/\, PO \,.\, BQ.$$

If the Euclidean length of OB is $r$, its hyperbolic length $\rho$ is given by

$$\rho = \tfrac{1}{2} \log (OBPQ)$$

$$= \tfrac{1}{2} \log \left| PB \,.\, OQ \,/\, PO \,.\, BQ \right|$$

$$= \tfrac{1}{2} \log [(a + r)\, a \,/\, a\, (a - r)]$$

$$= \tanh^{-1} (r/a),$$

assuming $R = 1$.

## 6. The Metric in Hyperbolic Polar Coordinates

Beltrami used the identity given in equation (7) to express the metric in hyperbolic polar coordinates $\rho$, $\phi$.

Given that

$$\rho = R \tanh^{-1}(r/a) ,$$

it follows that

$$d\rho = Ra\, dr / (a^2 - r^2)$$

or

$$\mathbf{a^2\, dr^2 / (a^2 - r^2)^2 = d\rho^2} \qquad \text{... (10)}$$

Furthermore from

$$r = a \tanh(\rho / R),$$

it is clear that

$$r^2 / (a^2 - r^2) = a^2 \tanh^2(\rho / R) / (a^2 - \tanh^2(\rho / R)$$

$$= \tanh^2(\rho / R) / \mathrm{sech}^2(\rho / R),$$

so

$$\mathbf{r^2 / (a^2 - r^2) = \sinh^2(\rho / R)} \qquad \text{... (11)}$$

Considering equation (3)

$$ds^2 = R^2 \{a^2\, dr^2 / (a^2 - r^2)^2 + r^2\, d\phi^2 / (a^2 - r^2) \}$$

in the light of (9) and (10), it is clear that

$$\mathbf{ds^2 = d\rho^2 + R^2 \sinh^2(\rho / R)\, d\phi^2,} \qquad \text{... (12)}$$

which expresses the metric in hyperbolic polar coordinates.

## 7. The Circumference of a Circle

Beltrami also provided the formula for the circumference $C$ of a circle on the disc,

$$\boldsymbol{C = 2\pi R \sinh(\rho/R)}, \qquad \textbf{... (13)}$$

where $\rho$ is the hyperbolic radius.

This identity, one of the most important formulas in hyperbolic geometry, can be derived from the metric given in equation (3):

$$ds^2 = R^2\{a^2\, dr^2 / (a^2 - r^2)^2 + r^2\, d\phi^2 / (a^2 - r^2)\}.$$

In a circle of radius $r$, $r$ is constant, so $dr = 0$, implying that

$$ds^2 = R^2 r^2\, d\theta^2 / (a^2 - r^2).$$

or

$$ds = R\, r\, d\theta / (a^2 - r^2)^{1/2}.$$

Integration then gives

$$C = R \int_0^{2\pi} r\, d\theta / (a^2 - r^2)^{1/2}$$

$$= 2\pi r R / (a^2 - r^2)^{1/2}.$$

But $r = a \tanh(\rho/R)$, whence

$$C = 2\pi Ra \tanh(\rho/R) / [a^2 - a^2 \tanh^2(\rho/R)]^{1/2}$$

$$= 2\pi R \tanh(\rho/R)\, \mathrm{sech}(\rho\ R)$$

$$= 2\pi R \sinh(\rho/R).$$

This identity constitutes the first of Beltrami's proofs that figures on the disc exhibit the properties of the hyperbolic plane.

## 8. The relationship between Angles on the Pseudosphere and Angles on the Disc

Next Beltrami obtained the relationship between the angle between a pair of geodesics on the pseudosphere and the angle between the chords which represent these geodesics on the disc.

Let $(u, v)$ be a point on the pseudosphere and let $(U, V)$ be a point on a geodesic emanating from it. Let the equations of two such geodesics be

$$V - v = n\,(U - u)$$

and

$$V - v = m\,(U - u).$$

At the point $(u, v)$ the geodesics follow directions given by

$$n = dv_1/du_1$$

and

$$m = dv_2/du_2\,,$$

so that $dv_1 = n\,du_1$ and $dv_2 = m\,du_2$.

The tangent vectors representing these directions are

$$d\mathbf{r}_1 = \mathbf{r_u}\,du_1 + \mathbf{r_v}\,dv_1$$

and

$$d\mathbf{r}_2 = \mathbf{r_u}\,du_2 + \mathbf{r_v}\,dv_2.$$

If $\alpha$ is the angle between these tangent vectors, then

$$\tan\alpha = (dr_1\ dr_2 \sin\alpha) / (dr_1\,dr_2 \cos\alpha)$$

$$= |d\mathbf{r}_1 \times d\mathbf{r}_2| / (d\mathbf{r}_1 \,.\, d\mathbf{r}_2)$$

$$= |(\mathbf{r_u}\,du_1 + \mathbf{r_v}\,dv_1) \times (\mathbf{r_u}\,du_2 + \mathbf{r_v}\,dv_2)| / [\,(\mathbf{r_u}\,du_1 + \mathbf{r_v}\,dv_1) \,.\, (\mathbf{r_u}\,du_2 + \mathbf{r_v}\,dv_2)\,].$$

But $\mathbf{r_u} \times \mathbf{r_u} = \mathbf{r_v} \times \mathbf{r_v} = 0$ and $(\mathbf{r}_u \times \mathbf{r}_u) = -(\mathbf{r}_v \times \mathbf{r}_u)$, so the numerator of the previous equation is

$$(\mathbf{r_v}\,dv_1 \times \mathbf{r_u}\,du_2) - (\mathbf{r_v}\,dv_2 \times \mathbf{r_u}\,du_1) = (dv_1\,du_2 - dv_2\,du_1)\,(\mathbf{r}_v \times \mathbf{r}_u)$$

$$= (n\,du_1\,du_2 - m\,du_1\,du_2)\,(\mathbf{r}_v \times \mathbf{r}_u)$$

$$= (n - m)\,du_1\,du_2\,(\mathbf{r}_v \times \mathbf{r}_u),$$

given that $n = dv_1/du_1$ and $dv_2 = m\,du_2$.

The denominator is

$$(\mathbf{r_u}\, du_1 + \mathbf{r_v}\, dv_1) \,.\, (\mathbf{r_u}\, du_2 + \mathbf{r_v}\, dv_2)$$

$$= \mathbf{r}_u \,.\, \mathbf{r}_u\; du_1\, du_2 + \mathbf{r}_v \,.\, \mathbf{r}_u\, n\, du_1\, du_2 \;+ \mathbf{r}_u \,.\, \mathbf{r}_v\, m\, du_1\, du_2 \;+ \mathbf{r}_v \,.\, \mathbf{r}_v\, nm\; du_1\, du_2 \,.$$

However $\mathbf{r}_u \,.\, \mathbf{r}_u = E$, $\mathbf{r}_u \,.\, \mathbf{r}_v = F$, $\mathbf{r}_v \,.\, \mathbf{r}_v = G$ and $\mathbf{r}_v \times \mathbf{r}_u = (EG - F^2)^{1/2}$, so the denominator is

$$[E + (n + m)\, F + nm\, G]\, du_1\, du_2\,.$$

Hence, with $du_1\, du_2$ cancelling out,

$$tan\,\alpha = (n - m)\; (EG - F^2)^{1/2} \,/\, [E + (n + m)\, F + nm\, G],$$

which is the formula quoted by Beltrami.

Into this formula he substituted the actual values of $E$, $F$ and $G$

$$E = [R^2/\, (a^2 - u^2 - v^2)^2]\; (a^2 - v^2)$$

$$F = [R^2/\, (a^2 - u^2 - v^2)^2]\; uv$$

$$G = R^2/\, (a^2 - u^2 - v^2)^2]\; (a^2 - u^2)$$

to obtain

$$\boldsymbol{tan\,\alpha = a\, (n - m)\, w \,/\, [\, (1 + mn)\, a^2 - (v - mu)\, (v - nu)\, ]} \qquad \textbf{... (15)}$$

where $w = (a^2 - u^2 - v^2)^{1/2}$.

Finally Beltrami considered the angle $\alpha'$ between the chords on the disc which maps the angle $\alpha$ between the corresponding geodesics on the pseudosphere. Letting $\mu$ and $\nu$ denote the angles which the chords make with the $x$-axis of the disc, he obtained

$$m = tan\,\mu$$

and

$$n = tan\,\nu,$$

where $m$, $n$ are the gradients of the two chords and $\alpha' = \nu - \mu$. From equation (15) he derived the identity relating the angle $\alpha$ to the angle $\alpha'$:

$$\boldsymbol{tan\,\alpha = a\, w\, sin\,\alpha' \,/\, [\, a^2 cos\,\alpha' - (v\, cos\,\mu - u\, sin\;\mu\,)\, (v\, cos\;\nu\; - u\, sin\;\nu\,)\, ]} \qquad \textbf{... (16)}$$

This equation can be justified as follows. Multiply the numerator of equation (15) by $\cos\nu \cos\mu$, to obtain

$$aw\,[\tan\nu \cos\nu \cos\mu - \tan\mu \cos\nu \cos\mu] = aw\,[\sin\nu\ \cos\mu - \cos\nu \sin\mu]$$

$$= aw\,[\sin\nu\ \cos\mu - \cos\nu \sin\mu]$$

$$= aw \sin(\nu - \mu)$$

$$= aw \sin\alpha',$$

which is the numerator of equation (16).

Multiply the denominator of equation (15) by the same factor $\cos\nu \cos\mu$, to obtain

$$\cos\nu \cos\mu\, [\,(1 + mn)\, a^2 - (v - mu)\,(v - nu)\,]$$

$$= \cos\nu \cos\mu\, [\,(1 + \tan\mu\ \tan\nu)\, a^2 - (v - u \tan\mu)\,(v - u \tan\nu)\,]$$

$$= [\, a^2 (\cos\nu \cos\mu + \sin\nu \sin\mu) - v^2 \cos\nu \cos\mu + uv \cos\nu \sin\mu + uv \cos\mu \sin\nu - u^2 \sin\mu \sin\nu$$

$$= a^2 \cos(v - \mu) - (v \cos\mu - u \sin\mu)\,(v \cos\nu - u \sin\nu)$$

$$= a^2 \cos\alpha' - (v \cos\mu - u \sin\mu)\,(v \cos\nu - u \sin\nu).$$

which is the denominator of equation (16), as required.

Beltrami noted that this expression for $\tan a$ in equation (16) is non-zero except when $w = (a^2 - u^2 - v^2)^{1/2} = 0$, that is to say, when the chords meet at the circumference, in which case the original geodesics only meet at an infinitely distant point.

He concluded that

1) Two distinct chords which intersect within the disc correspond to geodesics on the pseudosphere which meet at a finite point at angle $0 < a < 180°$

2) Two distinct chords which meet at the circumference of the disc correspond to geodesics which converge to a single infinite point and are inclined at an angle of zero

3) Two distinct chords which intersect outside the disc (i.e. which are parallel) correspond to geodesics which have no common point on the surface.

## 9. The Angle of Parallelism

At this point in his essay Beltrami sought concrete proofs that pseudospherical geometry agrees with hyperbolic planimetry. To this end he used chords on the disc, representing geodesics on the pseudosphere, to derive the hyperbolic angle of parallelism

$$\boldsymbol{\cot\theta = \sinh(\rho/R)} \qquad \textbf{... (17)}$$

In the figure below, $r$ is the Euclidean length of segment OB which is perpendicular to chord ST, while $\theta$ is the magnitude of angle SOB at the centre of the unit disc. The radius OS meets the chord ST at its endpoint S, regarded as infinitely far from the centre in terms of hyperbolic distance. OS is defined as parallel to ST since it separates those lines emanating from O such as OQ which intersect ST from those like OP which do not. Angle SOB = $\theta$ is the angle of parallelism corresponding to the segment OB:

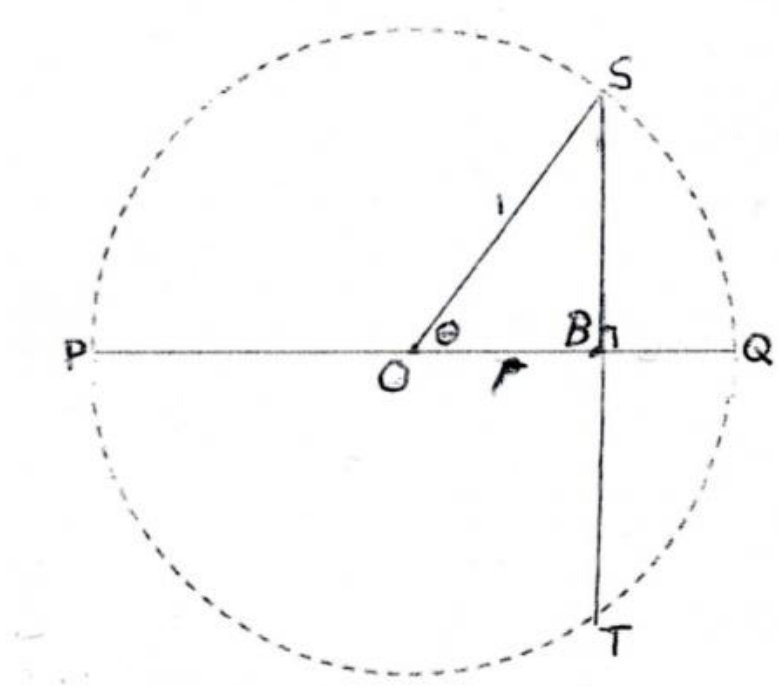

By equation (8), assuming that $a = 1$, it is clear that OB = $r = \tanh(\rho/R)$, which implies that

$$\text{SB} = [1 - \tanh^2(\rho/R)]^{1/2}$$
$$= \text{sech}(\rho/R).$$

Then

$$\tan\theta = \text{SB} / \text{OB}$$
$$= \text{sech}(\rho/R) / \tanh(\rho/R)$$
$$= \text{cosech}(\rho/R)$$

or

$$\cot\theta = \sinh(\rho/R).$$

Beltrami also showed that this identity is equivalent to

$$\tan \tfrac{1}{2}\theta = e^{-\rho/R},$$

the form in which it appears in Lobachevsky's work.

### 10. The Area of a Triangle

Beltrami also showed that the angle sum of a triangle on the disc is less than two right angles and that its area is proportional to its defect, two of the most important theorems of hyperbolic geometry.

He postulated a right angled (isosceles) triangle on where the base lies along the fundamental geodesic $v = 0$, the height is the perpendicular $u$ = constant and the hypotenuse issuing from

$$v = u \tan \mu .$$

The other acute angle of this triangle is denoted by $\mu'$. Suppose that this triangle is mapped on to Beltrami disc:

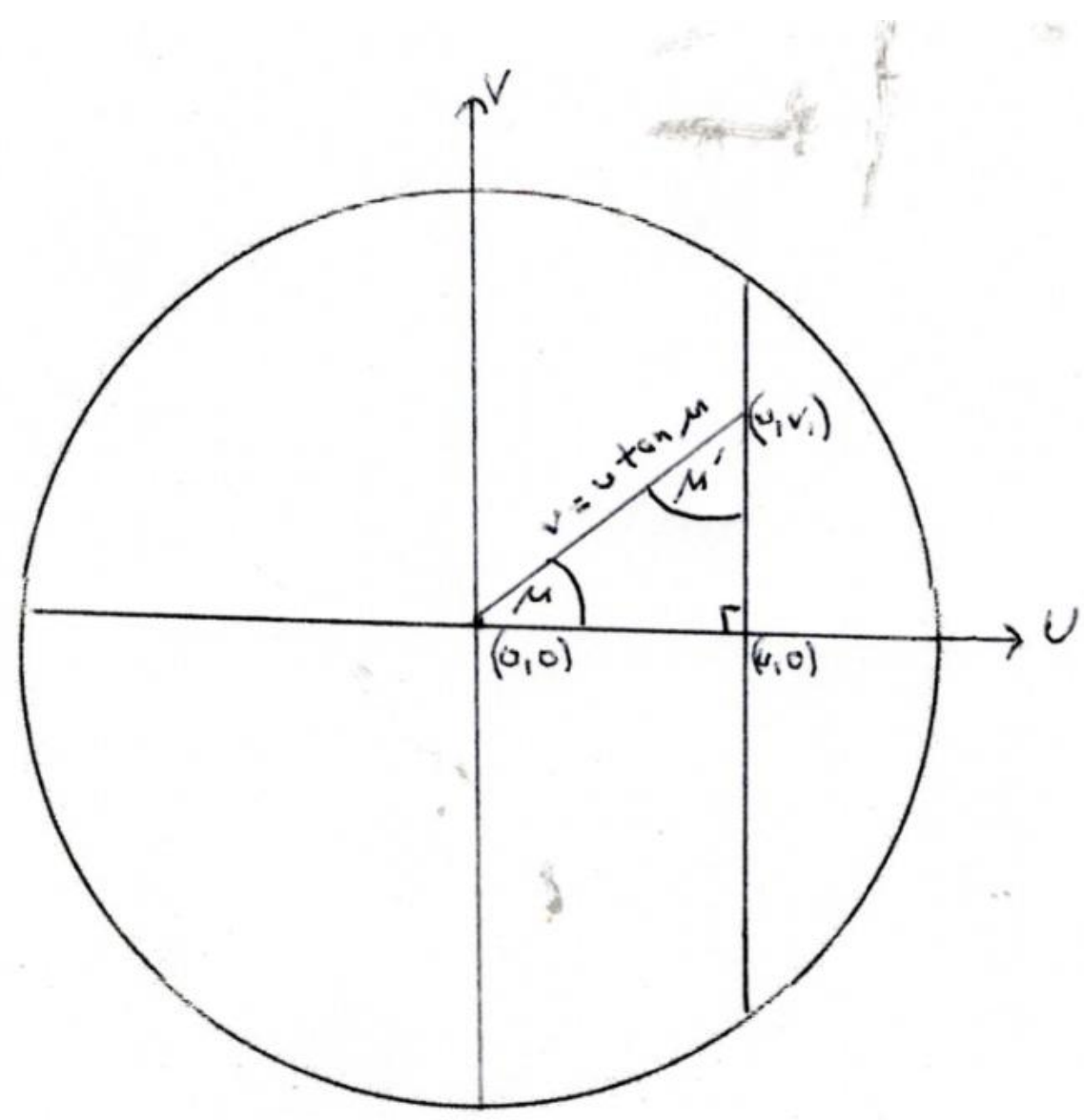

Viewed in Euclidean terms the acute angles of the triangle are respectively $\mu$ at the origin and $90° - \mu$ at the circumference. But the Beltrami map is not conformal except at the origin. So while in hyperbolic terms the acute angle at the centre is $\mu$, it cannot be assumed that the other acute angle $\mu'$ equals $90° - \mu$. To obtain the value of $\mu'$ it is necessary to use equation (16) which provides the hyperbolic angle if the Euclidean angle is given. However the equation needs to be slightly modified since the Euclidean angle in this case is not $\mu$ but the complement of $\mu$. The result is that $\sin \alpha'$ in the formula becomes $\cos \mu$ and $\cos \alpha'$ becomes $\sin \mu$ . Hence

$$\tan \mu' = a\, w \cos \mu \;/\; [\, a^2 \sin \mu - (v \cos \mu - u \sin \mu)\,(v \cos \nu - u \sin \nu)\,]$$

$$= a\, (a^2 - u^2 - v^2)^{1/2} \cos \mu \;/\; a^2 \sin \mu,$$

so

$$\boldsymbol{\tan \mu' = (a^2 \cos^2 \mu - u^2)^{1/2} \,/\, a \sin \mu} \qquad \textbf{... (18)}$$

The second line of this derivation, in which the factors ($v \cos \nu - u \sin \nu$) have been discarded, follows from the fact that the arms of the angle $\mu'$ are inclined to the x-axis at angles $\mu$ and $\nu$ where $\nu$ = R and $\cos \nu = 0$. Furthermore $v = u \tan \mu$. Therefore

$$(v \cos \mu - u \sin \mu)(v \cos \nu - u \sin \nu) = -u \sin \nu (u \tan \mu \cos \mu - u \sin \mu)$$

$$= -uv \sin \nu \sin \mu + uv \sin \nu \sin \mu$$

$$= 0.$$

The last line also follows from the fact that $v = u \tan \mu$, whence

$$(a^2 - u^2 - v^2)^{1/2} \cos \mu = [(a^2 - u^2 - v^2) \cos^2 \mu]^{1/2}$$

$$= [a^2 \cos^2 \mu - u^2 (\cos^2 \mu + \sin^2 \mu)]^{1/2}$$

$$= (a^2 \cos^2 \mu - u^2)^{1/2}.$$

It then follows by differentiating both sides of equation (18) that

$$\mathbf{d\mu' = -[a \sin \mu\, u\, du] / (a^2 - u^2)(a^2 \cos^2 \mu - u^2)^{1/2}} \quad \textbf{... (19)}$$

which, as Beltrami noted, expresses the increment received by angle $\mu'$ when the side of the triangle opposite angle $\mu$ moves parallel to itself towards the origin (angle $\mu$ being held constant.

At this point Beltrami considered the (hyperbolic) element of area within the triangle given by

$$du\, dv\, (EG - F^2)^{1/2} = R^2 a\, du\, dv / (a^2 - u^2 - v^2)^{3/2}.$$

First he integrated the element of area with respect to $v$ from $v = 0$ to $v = u \tan \mu$ (i.e. vertically), keeping $u$ constant, obtaining

$$R^2 [a \sin \mu\, u\, du] / (a^2 - u^2)(a^2 \cos^2 \mu - u^2)^{1/2}.$$

This expression can be justified as follows. Writing $R^2 a\, du$ as M and ($a^2 - u^2$) as $k$, the integral

$$J = \int R^2 a\, du\, dv / (a^2 - u^2 - v^2)^{3/2}$$

takes the form

$$J = M \int dv / (k^2 - v^2)^{3/2}.$$

Let $v = k \sin \theta$. Then $dv = k \cos \theta\, d\theta$ and $(k^2 - v^2)^{3/2} = (k^2 - k^2 \sin^2 \theta)^{3/2} = k^3 \cos^3 \theta$. It follows that

$$J = M \int k \cos \theta\, d\theta \;/\; k^3 \cos^3 \theta$$

$$= M \int d\theta \;/\; k^2 \cos^2 \theta$$

$$= M \int \sec^2 \theta\, d\theta / k^2$$

$$= M \tan \theta / k^2.$$

But $tan\ \theta = v / (k^2 - v^2)^{1/2}$, so

$$J = (R^2\ a\ du)\ v / (a^2 - u^2)\ (a^2 - u^2 - v^2)^{1/2}\ \Big|_0^{u\ tan\ \mu}$$

Now when $v = 0$, the value of the integral is zero, so it is only necessary to substitute $u\ tan\ \mu$ for $v$, yielding

$$J = (R^2\ a\ du)\ u\ tan\ \mu / (a^2 - u^2)\ (a^2 - u^2 - u^2\ tan^2\ \mu)^{1/2}$$

$$= (R^2\ a\ du)\ u\ tan\ \mu / (a^2 - u^2)\ [a^2 - u^2\ (1 + tan^2\ \mu)]^{1/2}$$

$$= (R^2\ a\ du)\ u\ tan\ \mu / (a^2 - u^2)\ [a^2 - u^2\ sec^2\ \mu]^{1/2}$$

$$= (R^2\ a\ du)\ u\ sin\ \mu / (a^2 - u^2)\ [a^2\ cos^2\ \mu - u^2]^{1/2},$$

as required. But by the miracle of equation (18)

$$\mathbf{(R^2\ a\ du)\ u\ sin\ \mu / (a^2 - u^2)\ [a^2\ cos^2\ \mu - u^2]^{1/2} = -R^2\ d\mu'} \qquad \text{... (20)}$$

Beltrami then turned his attention to integrating the element of area with respect to $u$ (i.e. horizontally). But rather than integrating the left hand side of the previous equation, he integrated the right hand side from $\mu' = 90° - \mu$ to $\mu' = \mu$ . When angle $\mu'$ equals $90° - \mu$, the angle at the centre being $\mu$, the triangle is Euclidean, which in the context of hyperbolic geometry is only possible for infinitely small triangles i.e. when $u = 0$. Beltrami did not explain the significance of the other limit of his integral. But if the triangle is assumed to be isosceles, the value $\mu' = \mu$ guarantees that the sides about the right angle are equal. (This is true because Euclid's isosceles angle theorem is valid in the hyperbolic plane). So when $\mu' = \mu$, $u = u_1$. Thus in integrating the right hand side of (15) with respect to $\mu'$, Beltrami effectively integrated the left hand side with respect to $u$ all along the base of the triangle, obtaining

$$-R^2 \int d\mu' = -R^2 \mu'\ \Big|_{90° - \mu}^{\mu}$$

$$= R^2\ (90° - \mu - \mu).$$

But at the maximum value of $u$, the angle $\mu' = \mu$, so one of the $\mu$ s can be written as $\mu'$. Thus Beltrami finally obtained the area of the triangle as

$$R^2\ (\pi/2 - \mu - \mu').$$

For any positive area, this expression implies firstly that the sum of the angles of a right angled triangle is less than $\pi$ and secondly that the area of the triangle is proportional to its defect. At this point Beltrami noted that, since any triangle can divided into two right angled triangles, it follows that the area of any triangle ABC is given by

$$\mathbf{R^2\ [\pi - (\angle A + \angle B + \angle C)]} \qquad \text{... (21)}$$

Beltrami also noted that in the case of a triangle whose vertices lie on the circumference, that is to say, where the sides are asymptotically parallel and the sum of the angles is zero, the area by equation (16) will be $R^2\pi$. This is the triangle of maximum size of hyperbolic geometry, a figure whose area was evaluated by Bolyai as $\pi k^2$, where $k$ is the linear constant of the hyperbolic plane, equivalent to $R$ , the reciprocal (with its sign changed) of the curvature of the surface.

In fact Beltrami only proved his theorem on the angle sum of a triangle in respect of an isosceles right angled triangle. But he could have invoked Saccheri's theorem that if in a single case the sum of the angles of a triangle is less than two right angles, then it is true in every case.

## 11. Geodesic Circles

A geodesic circle is defined as the set of points lying at a constant geodetic distance $\rho$ from a given point. Beltrami noted that a geodesic circle centred at an arbitrary point $(u_0, v_0)$ with hyperbolic radius $\rho$ will satisfy

$$\boldsymbol{(a^2 - uu_o - vv_0) / (a^2 - u^2 - v^2)^{1/2}(a^2 - u_o^2 - v_0^2)^{1/2} = \cosh(\rho / R)} \qquad \textbf{... (22)}$$

He derived this equation in a supplementary Note, but it can be seen to be true in the case where $(u_0, v_0)$ is (0, 0):

$$(a^2 - uu_o - vv_0) / (a^2 - u^2 - v^2)^{1/2} (a^2 - uu_o - vv_0)^{1/2} = a^2 / a (a^2 - u^2 - v^2)^{1/2}$$

$$= a / (a^2 - r^2)^{1/2}$$

$$= \cosh(\rho / R),$$

by equation (9).

In order to derive the equation of a circle on the disc, Beltrami first obtained the equation of a system of chords emanating from an arbitrary point $(u_0, v_0)$ and then obtained the equation of an orthogonal system, since, clearly, the circumference of any circle will be orthogonal to the radii emanating from the centre.

He gave

$$v - v_0 = k (u - u_0)$$

as the equation of chords emanating from the point $(u_0, v_0)$, which is real or ideal depending on whether $u_0^2 + v_0^2$ is less than or greater than $a^2$. It follows that

$$k = (v - v_0) / (u - u_0)$$

and given that $dv = k\, du$, it is clear that the differential equation of the system is

$$du / (u - u_0) = dv / (v - v_0).$$

This implies that the direction vectors of $u$, $v$ are respectively proportional to $(u - u_0)$, $(v - v_0)$. The differential equation of the orthogonal system with direction vectors $du$, $dv$ is then by virtue of a well known theorem given by

$$\boldsymbol{E (u - u_0)\, du + F (v - v_0)\, du + F (u - u_0)\, dv + G (v - v_0)\, dv = 0} \qquad \textbf{... (23)}$$

Given the actual values of $E$, $F$ and $G$, Beltrami stated that this equation is equivalent to

$$\boldsymbol{d \{ (a^2 - uu_o - vv_0) / (a^2 - u^2 - v^2)^{1/2} \} = 0} \qquad \textbf{... (24)}$$

whence by integration

$$\boldsymbol{(a^2 - uu_o - vv_0) / (a^2 - u^2 - v^2)^{1/2} = C} \qquad \textbf{... (25)}$$

The equivalence of equations (23) and (24) can be proved as follows. The values of $E$, $F$ $G$ are

$$E = R^2 / (a^2 - v^2) / (a^2 - u^2 - v^2)^2$$

$$F = R^2\, uv / (a^2 - u^2 - v^2)^2$$

$$G = R^2 / (a^2 - u^2) / (a^2 - u^2 - v^2)^2.$$

Therefore the left hand side of equation (22) becomes

$$R^2 / (a^2 - u^2 - v^2)^2 \{ ( a^2 - v^2) (u - u_0)\, du + uv (v - v_0)\, du + uv (u - u_0)\, dv + (a^2 - u^2) (v - v_0)\, dv$$

or

$$R^2 / (a^2 - u^2 - v^2)^2 \{ [a^2 u - a^2 u_0 + v^2 u_0 - uvv_0]\, du + [a^2 v - a^2 v_0 + u^2 u_0 - uvu_0]\, dv \}.$$

At the same time the total differential $df$ of $(a^2 - uu_o - vv_0) / (a^2 - u^2 - v^2)^{1/2}$ in equation (23) is

$$df = \partial f / \partial u \; du + \partial f / \partial v \; dv$$

$$= (1 / (a^2 - u^2 - v^2)^{3/2} \{ [a^2 u - a^2 u_0 + v^2 u_0 - uvv_0]\, du + [a^2 v - a^2 v_0 + u^2 u_0 - uvu_0]\, dv \}.$$

Since the right hand side in both equations is zero, while the expressions in curly brackets are identical, the two equations are equivalent, as stated.

Comparing equations (22) and (25) Beltrami obtained

$$\boldsymbol{\cosh (\rho / R) = C / (a^2 - u^2 - v^2)^{1/2}} \qquad \textbf{... (26)}$$

He noted that in the case of a circle with a real centre and a real radius $\rho$, $C$ can never be zero. This is because *cosh x* cannot equal zero for real $x$.

**12. Equidistants**

Beltrami pointed out that when the centre of the circumference is an ultra-ideal point (i.e. outside the disc), the constant *C* can take the value zero, for then equation (25) becomes

$$(a^2 - uu_o - vv_0) = 0 \qquad \dots (27)$$

which, by a well known theorem, represents the polar of the external point ($u_0$, $v_0$). In this case the pole at ($u_0$, $v_0$) is the ultra-ideal centre of the curve equidistant at every point to the polar, an equidistant curve being a key feature of hyperbolic geometry.

Beltrami chose as his polar the diameter $v = 0$ (whose pole would lie at an infinitely distant point outside the disc) and drew the perpendicular from Q ($u$, 0) on the diameter, meeting the equidistant curve at P ($u$, $v$). The endpoints of the chord on which PQ lies are A and B:

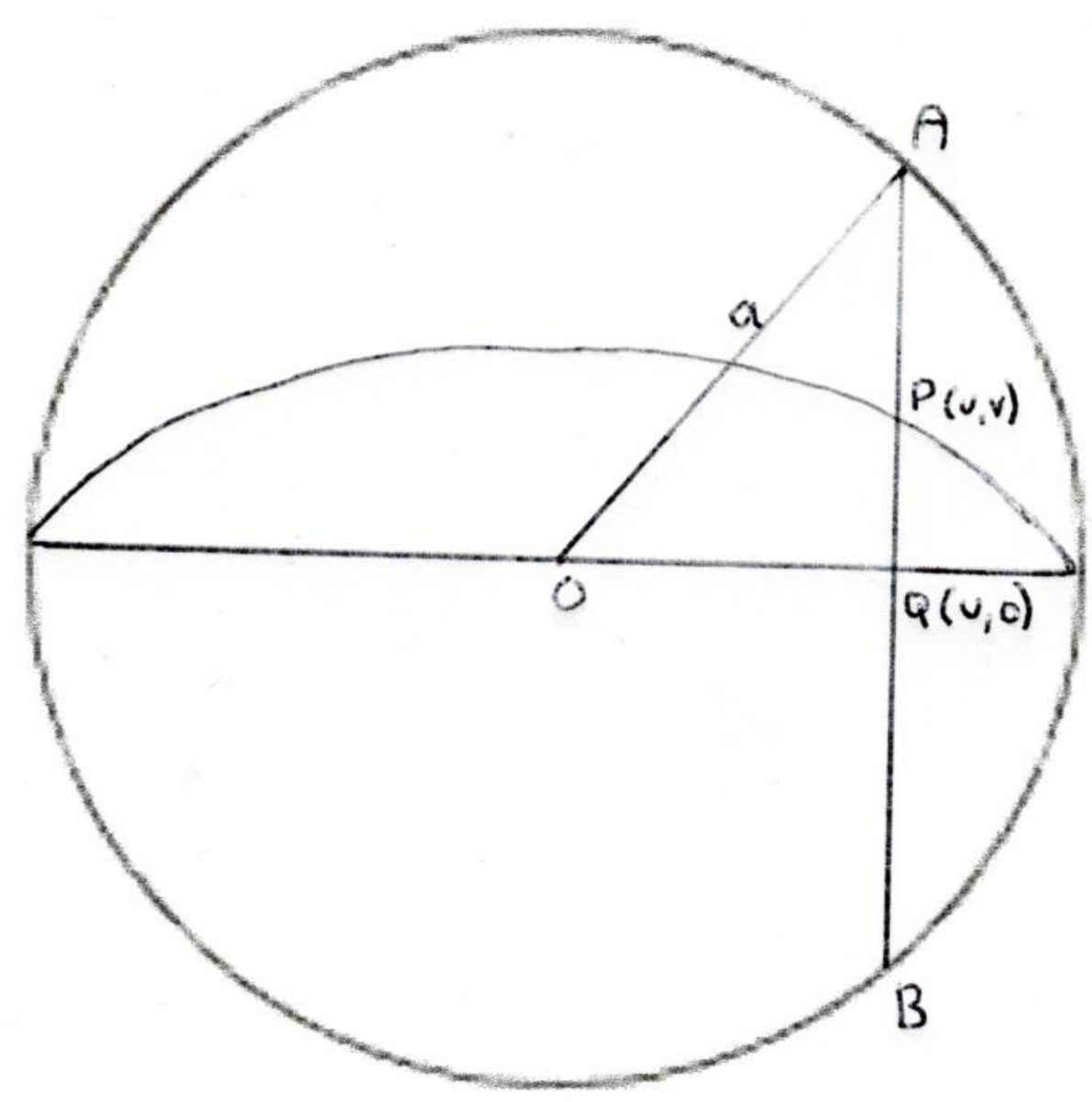

He denoted the hyperbolic length of PQ by $\xi$, where by section 5 above

$$\xi = R \tfrac{1}{2} \log (\mathrm{PQAB})$$

$$= R \tfrac{1}{2} \log \left| \mathrm{BP} \,.\, \mathrm{QA} / \mathrm{BQ} \,.\, \mathrm{PA} \right|$$

$$= R \tfrac{1}{2} \log \left| [ (a^2 - u^2)^{1/2} + v] / [ (a^2 - u^2)^{1/2} - v] \right| .$$

If $\xi$ is supposed constant, then

$$\xi = R \tfrac{1}{2} \log \left| [ (a^2 - u^2)^{1/2} + v] / [ (a^2 - u^2)^{1/2} - v] \right| \qquad \dots (28)$$

is the equation of any equidistant curve having its ultra-ideal centre at ($u_0$, $v_0$).

He also denoted the hyperbolic length of OQ by $\eta$, where by the same theorem

$$\eta = R \tfrac{1}{2} \log \left| (a + u) / (a - u) \right| \qquad \dots (29)$$

By equation (8)

$$\boldsymbol{u = a \tanh (\eta/R)} \qquad \text{... (30)}$$

By the same token

$$v = (a^2 - u^2)^{1/2} \tanh (\xi/R).$$

But by equation (9)

$$(a^2 - u^2)^{1/2} = a / \cosh (\eta/R) ,$$

whence

$$\boldsymbol{v = a \tanh (\xi/R) / \cosh (\eta/R)} \qquad \text{... (31)}$$

When the values of $u$ and $v$ are substituted into the $w^2 = (a^2 - u^2 - v^2)$, then

$$w^2 = a^2 - a^2 \tanh^2 (\eta/R) - a^2 \tanh^2 (\xi/R) / \cosh^2 (\eta/R)$$

$$= a^2 (1 - \tanh^2 (\eta/R) ) - a^2 \tanh^2 (\xi/R) / \cosh^2 (\eta/R)$$

$$= a^2 / \cosh^2 (\eta/R) - a^2 \tanh^2 (\xi/R) / \cosh^2 (\eta/R)$$

$$= a^2 / \cosh^2 (\eta/R) [ 1 - \tanh^2 (\xi/R) ]$$

$$= a^2 / (\cosh^2 (\eta/R) (\cosh^2 (\xi/R).$$

Substituting these values of $u$, $v$ and $w$ into the metric given in equation (1), Beltrami obtained

$$\boldsymbol{ds^2 = d\xi^2 + \cosh^2 (\xi/R) d\eta^2} \qquad \text{... (32)}$$

which he identified as the metric of a surface of revolution.

### 13. Horocycles

Having discussed circumferences (to use Beltrami's generic term) with a real centre inside the disc and circumferences with an ultra-ideal centre outside the disc, Beltrami considered the intermediate circumferences with an ideal centre lying on the boundary of the disc. These take the form of nested circles (or ellipses) which are tangent from the inside to a point on the boundary of the disc. He showed that these circles map concentric horocycles on the pseudospherical surface (such as parallels of latitude on a pseudosphere) while the chords emanating from the point of tangency and orthogonal to the circles map the axes of the horocycles (such as the meridians on a pseudosphere):

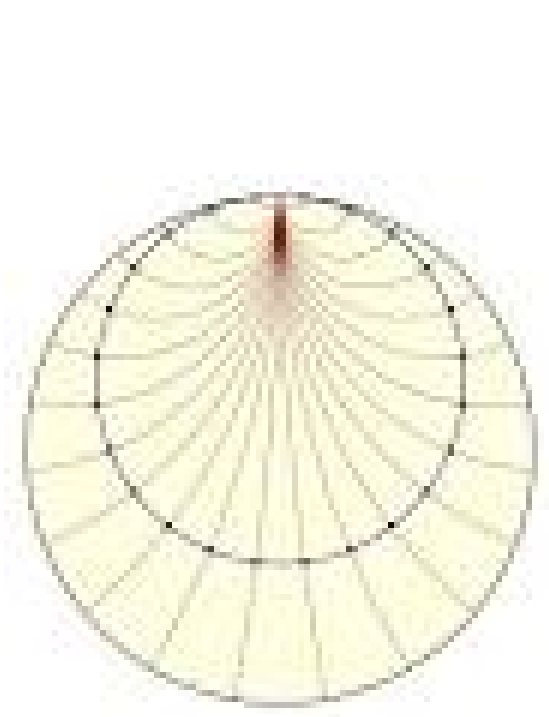

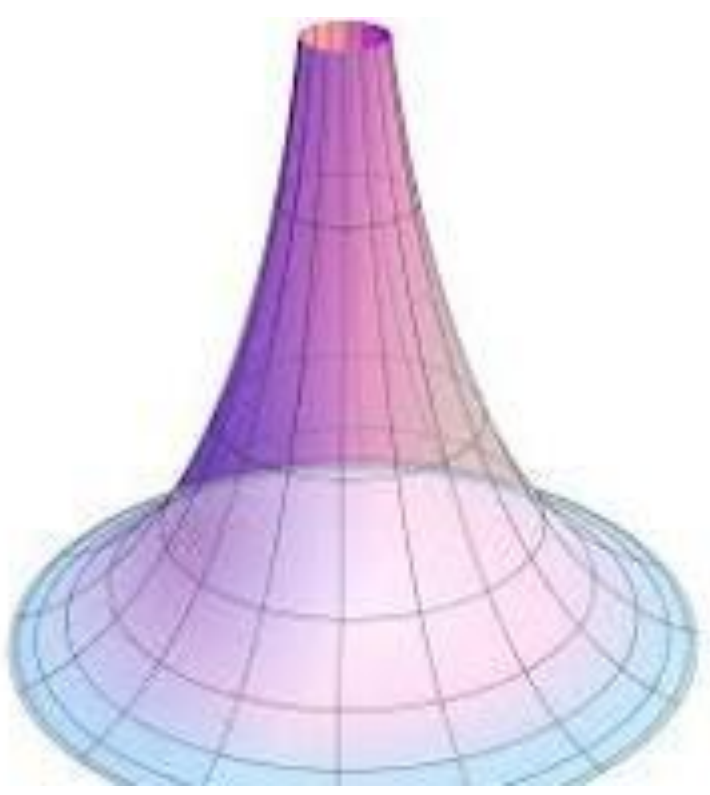

His starting point was equation (25)

$$(a^2 - uu_o - vv_0) / (a^2 - u^2 - v^2)^{1/2}(a^2 - u_o^2 - v_0^2)^{1/2} = \cosh(\rho/R),$$

or

$$\mathbf{(a^2 - uu_o - vv_0) / (a^2 - u^2 - v^2)^{1/2} = w_0 \cosh(\rho/R),} \qquad \textbf{... (33)}$$

where $w_0 = (a^2 - u_o^2 - v_0^2)^{1/2}$.

When the centre of the circle ($u_0$, $v_0$) moves to the boundary of the disc, $w_0$ converges to zero. At the same time the hyperbolic radius $\rho$ increases without limit. If now $\rho$ is written as $\rho' - \rho$, where $\rho$ remains finite while $\rho'$ increases indefinitely, the previous equation becomes

$$(a^2 - uu_o - vv_0) / (a^2 - u^2 - v^2)^{1/2} = w_0 \cosh[(\rho' - \rho)/R]$$

$$= (w_0/2)(e^{\rho'/R} e^{-\rho/R}) + (w_0/2)(e^{-\rho'/R} e^{\rho/R}).$$

Now the second term on the right hand side converges to zero since both $w^0$ and $e^{-\rho'/R}$ become infinitely small, while the first term remains finite (assuming that $w_0$ approaches zero as fast as $e^{\rho'/R}$ becomes infinite). In this case

$$\mathbf{(a^2 - uu_o - vv_0) / (a^2 - u^2 - v^2)^{1/2} = k\, e^{-\rho/R}} \qquad \textbf{... (34)}$$

where $k$ is a constant.

Beltrami used equation (34) and the following (mysterious) equation,

$$(u_0 v - uv_0) / (a^2 - uu_0 - vv_0) = \sigma / R, \qquad \text{... (35)}$$

to derive the well known formula for the linear element of the hyperbolic plane in horocyclic coordinates

$$ds^2 = d\rho^2 + e^{-2\rho/R} d\sigma^2. \qquad \text{... (36)}$$

Here $\rho$ denotes the perpendicular distance from a reference horocycle $\rho = 0$, while $\sigma$ denotes arc length along the reference horocycle, which is proportional to the (Euclidean) angle at the ideal point between the chords which delineate the arc.

Since a horocyclic arc remains at a constant distance from the ideal centre, the quantity $d\rho = 0$, implying that arc length along a horocycle is given by

$$ds = e^{-\rho/R} d\sigma.$$

In other words arc length along concentric horocycles diminishes by a factor $e^{-\rho/R}$ as the horocycles approach the boundary. So given two horocyclic arcs $s_1$, $s_2$ separated by a distance $\tau$, it follows that

$$s_2 = s_1 e^{-\tau/R}, \qquad \text{... (37)}$$

which accords with Lobachevsky's equation relating concentric horocyclic arcs on the hyperbolic plane.

**Appendix 1: Derivation of the Metric of a Surface of Negative Curvature**

. Suppose that the surface is a hyperboloid whose equation is

$$x^2 + y^2 - z^2 = - \mathrm{R}^2,$$

upon which a hyperbolic metric is imposed, namely

$$ds^2 = dx^2 + dy^2 - dz^2.$$

Suppose that a point ($x$, $y$, $z$) on the hyperboloid is projected along the ray passing through the origin until it intersects the plane the plane $z = R$. The intersection point ($u$, $v$) is obtained by writing

$$u = ax/z$$

$$v = ay/z,$$

for some positive $a$.

This ensures, since $x^2 + y^2 < z^2$ by virtue of the equation of the hyperboloid, that $u^2 + v^2 < a^2$. That is to say the point ($u, v$) will lie inside the Beltrami-Klein disc of radius $a$.

Now the definitions of $u$ and $v$ imply that

$$x = uz/a, y = vz/a$$

and so

$$x^2 = (uz/a)^2, y = (vz/a)^2.$$

These expressions for $x^2$ and $y^2$ are substituted into the equation of the hyperboloid to give

$$(uz/a)^2 + (vz/a)^2 - z^2 = - R^2,$$

whence

$$(z^2 /a^2)\,(u^2 + v^2 - a^2) = - R^2$$

or

$$z^2 = R^2a^2 / (a^2 - u^2 - v^2),$$

whence

$$z = Ra / (a^2 - u^2 - v^2)^{1/2}$$

or, writing $(a^2 - u^2 - v^2)^{1/2}$ as $w$,

$$z = Ra / w.$$

The differentials of $x$ , $y$ , $z$ are

$$dx = 1/a\ (z\ du + u\ dz)$$

$$dy = 1/a\ (z\ dv + v\ dz)$$

$$dz = Rau / w^3\ du + Rav / w^3\ dv$$

$$= [Ra / w^3]\ (u\ du + v\ dv).$$

However using the expressions for $z$ and $dz$,

$$dx = 1/a\ \ [Ra / w + uRa\ (u\ du + v\ dv) / w^3]$$

$$= R / w^3\ [\ (a^2 - v^2)\ du + uv\ dv].$$

Similarly

$$dy = R / w^3\ [\ (a^2 - u^2)\ dv + uv\ du].$$

Substituting $dx$, $dy$, $dz$ into the metric $ds^2 = dx^2 + dy^2 - dz^2$ gives

$$ds^2 = R^2\ [((a^2 - v^2)\ du + uv\ dv)^2 + ((a^2 - u^2)\ dv + uv\ du)^2 - a^2(u\ du + v\ dv)^2] / w^6$$

The coefficients of the terms in $du^2$, $dv^2$ and $du\ dv$ are

$$du^2: \quad (a^2 - v^2)^2 + u^2v^2 - a^2u^2 = (a^2 - v^2)\ w^2$$

$$dv^2: \quad (a^2 - u^2)^2 + u^2v^2 - a^2v^2 = (a^2 - u^2)\ w^2$$

$$du\ dv: \quad 2uv\ (a^2 - v^2) + 2uv\ (a^2 - u^2) - 2a^2uv = 2uv\ w^2.$$

Then

$$ds^2 = (R^2 / w^6)\ \ [\ (a^2 - v^2)\ w^2\ du^2 + 2uv\ du\ dv\ w^2 + (a^2 - u^2)\ w^2\ dv^2$$

$$= R^2\ [\ (a^2 - v^2)\ \ du^2 + 2uv\ du\ dv + (a^2 - u^2)\ \ dv^2] / w^4$$

or

$$ds^2 = R^2\ [\ (a^2 - v^2)\ \ du^2 + 2uv\ du\ dv + (a^2 - u^2)\ \ dv^2] / (a^2 - u^2 - v^2)^2\ ,$$

as required.

**Appendix 2: The Metric of a Hyperboloid of Negative Curvature in Polar Coordinates**

J. F. Barrett (2015) provides an alternative derivation of the metric of a hyperbolic surface in polar coordinates. He considers the upper sheet of a hyperboloid whose equation, given in Cartesian coordinates, is

$$x^2 + y^2 - z^2 = -R^2:$$

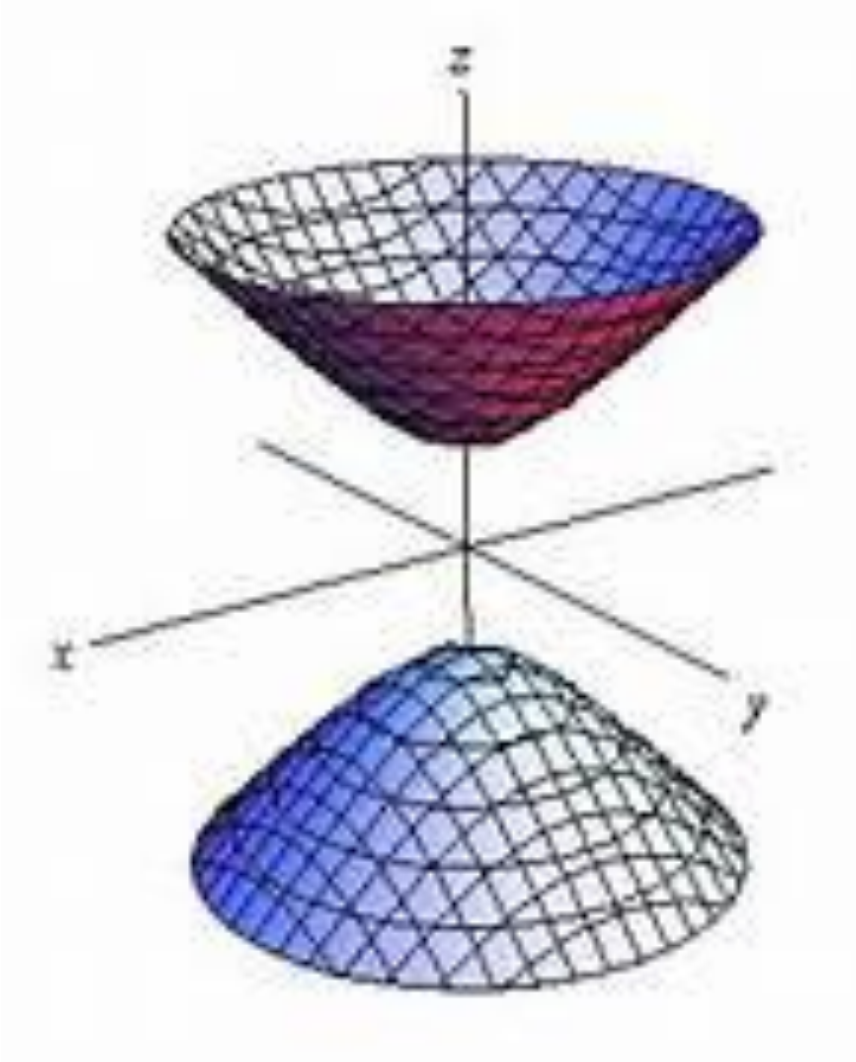

The *x*, *y*, *z* coordinates can be parametrised as follows:-[

$$x = R \sinh m \cos \theta$$

$$y = R \sinh m \sin \theta$$

$$z = R \cosh m,$$

where $m = p/R$, *p* being the length of the geodesic from the point (*x*, *y*, *z*) to the vertex, while θ is the azimuthal angle.

A hyperbolic metric arises if

$$ds^2 = dx^2 + dy^2 - dz^2.$$

Now

$$dx = R\,[\cosh m \cos \theta\, dm - \sinh m \sin \theta\, d\theta],$$

so

$$dx^2 = R^2\,[\cosh^2 m \cos^2 \theta\, dm^2 + \sinh^2 m \sin^2 \theta\, d\theta^2 \;- 2 \cosh m \cos \theta\; \sinh m \sin \theta\, dm\, d\theta],$$

while

$$dy = R\,[\cosh m \sin \theta\, dm + \sinh m \sin \cos d\theta].$$

Therefore

$$dy^2 = R^2 [cosh^2 \sin^2\theta\, dm^2 + sinh^2 m \cos^2\theta\, d\theta^2 + 2 \cosh m \sin\theta \sinh m \cos\theta\, dm\, d\theta,$$

and

$$dz^2 = R^2 sinh^2 \ dm^2.$$

Then, substituting into the metric,

$$ds^2 = R^2 [cosh^2 m \ dm^2 + sinh^2 m \ d\theta^2 - sinh^2 m\, dm^2]$$

or

$$\mathbf{ds^2 = R^2 [dm^2 + sinh^2 m\, d\theta^2]} \qquad ($$

a formula obtained by Beltrami.

**Appendix 3: Mapping Geodesics on a Surface to Straight Lines in the Plane**

In his first supplementary *Note* to his *Essay* Beltrami outlined the condition which a geodesic on a surface must satisfy if it is to be mapped to a straight line in the plane:

$$du \ d^2v - dv \ d^2u = 0.$$

This can be justified as follows. Supposing that $u$ and $v$ are functions of $t$, the previous differential equation is equivalent to

$$du/dt \,.\, d^2v/dt^2 - \mathrm{d}v/dt \,.\, \mathrm{d}^2u/dt^2 = 0$$

or

$$d/dt\, (\dot{u} / \dot{v}) = 0,$$

which implies that

$$\dot{u} / \dot{v} = C,$$

or

$$du/dv = C,$$

whence

$$u = Cv + \text{constant},$$

which is the equation of a straight line.

[1] There are however at least three mathematical typos in John Stillwell's excellent translation of Beltrami's *Essay* given in his *Sources of Hyperbolic Geometry*:

The first occurs on p.13 where the equation for $\tan \mu'$ should read

$$\tan \mu' = a\, w \cos \mu \;/\; [\, a^2 \sin \mu - (v \cos \mu - u \sin \mu)\, (v \cos \nu - u \sin \nu)\, ],$$

the translation omitting $u$ before $\sin \nu$ in the last bracket in the denominator.

The second occurs on p.17 where $u = a \tan \delta/R$ should be $u = a \tanh \delta/R$.

The third occurs in an equation found on p.19. The equation should read

$$d\mu' = -\,[a \sin \mu\, u\, du] \,/\, (a^2 - u^2)\, (a^2 \cos^2\mu - u^2)^{1/2},$$

with a minus sign at the beginning of the right hand side.